\documentclass[12pt]{amsart}
\usepackage{amsmath,amssymb}

\newtheorem{theorem}{Theorem}[section]
\newtheorem{lemma}[theorem]{Lemma}
\newtheorem{proposition}[theorem]{Proposition}
\newtheorem{corollary}[theorem]{Corollary}
\newtheorem{e-definition}[theorem]{Definition\rm}

\DeclareMathOperator{\LCD}{LCD}

\def \R {\mathbb{R}}

\def \Z {\mathbb{Z}}
\def \E {\mathbb{E}}
\def \P {\mathbb{P}}

\def \one {{\bf 1}}

\def \EE {\mathcal{E}}

\def \a {\alpha}

\def \g {\gamma}
\def \e {\varepsilon}

\def \d {\delta}

\def \l {\lambda}

\def \< {\langle}
\def \> {\rangle}
\def \^ {\widehat}

\def \dist {{\rm dist}}

\def \Span {{\rm span}}

\begin{document}

\title[The least singular value of a random matrix]
  {The least singular value of a random square matrix is $O(n^{-1/2})$}

\author{Mark Rudelson}
\author{Roman Vershynin}

\thanks{M.R. was supported by NSF DMS grant 0652684.
  R.V. was supported by the Alfred P.~Sloan Foundation
  and by NSF DMS grant 0652617.}

\address[M.R.]{Department of Mathematics,
   University of Missouri,
   Columbia, MO 65211, USA}
\email{rudelson@math.missouri.edu}   
   
\address[R.V.]{Department of Mathematics,
   University of California,
   Davis, CA 95616, USA}
\email{vershynin@math.ucdavis.edu}

\maketitle

\begin{abstract}
Let $A$ be a matrix whose entries are real i.i.d. centered random variables
with unit variance and suitable moment assumptions. Then the smallest singular value
$s_n(A)$ is of order $n^{-1/2}$ with high probability.
The lower estimate of this type was proved recently by the authors;
in this note we establish the matching upper estimate.
\end{abstract}

\section{Introduction}

Let $A$ be an $n \times n$ matrix whose entries are real i.i.d. centered random variables
with suitable moment assumptions.
Random matrix theory studies the distribution of the {\em singular values} $s_k(A)$,
which are the eigenvalues of $|A| = \sqrt{A^*A}$ arranged in the non-increasing order.
In this paper we study the magnitude of the smallest singular value $s_n(A)$,
which can also be viewed as the reciprocal of the spectral norm:
\begin{equation}                                \label{sn}
  \qquad s_n(A) = \inf_{x:\; \|x\|_2 = 1} \|Ax\|_2
  = 1/\|A^{-1}\|.
\end{equation}
Motivated by numerical inversion of large matrices,
von~Neumann and his associates speculated that
\begin{equation}                \label{sim}
  s_n(A) \sim n^{-1/2} \quad \text{with high probability.}
\end{equation}
(See \cite{vN}, pp.~14, 477, 555). A more precise form of this
estimate was conjectured by Smale and proved by Edelman \cite{E 88}
for Gaussian matrices $A$.
For general matrices, conjecture \eqref{sim} had remained open
until we proved in \cite{RV square} the lower bound
$s_n(A) = \Omega(n^{-1/2})$. In the present paper, we shall prove
the corresponding upper bound $s_n(A) = O(n^{-1/2})$, thereby
completing the proof of \eqref{sim}.

\begin{theorem}[Fourth moment]       \label{t: 4}
  Let $A$ be an $n \times n$ matrix whose entries are i.i.d.
  centered random variables with unit variance and fourth moment
  bounded by $B$.
  Then, for every $\d > 0$ there exist $K > 0$ and $n_0$
  which depend (polynomially) only on $\d$ and $B$, and such that
  $$
  \P \big(s_n(A) > K n^{-1/2} \big) \le \d
  \qquad \text{for all $n \ge n_0$}.
  $$
\end{theorem}

\noindent {\bf Remark.}
  The same result but with the reverse estimate,
  $\P \big(s_n(A) < K n^{-1/2} \big) \le \d$,
  was proved in \cite{RV square}. Together, these two estimates
  amount to \eqref{sim}.
  
\medskip

Under more restrictive (but still quite general) moment assumptions,
Theorem~\ref{t: 4} takes the following sharper form.
Recall that a random variable $\xi$ is called {\em subgaussian} if its tail is
dominated by that of the standard normal random variable:
there exists $B > 0$ such that
$\P (|\xi| > t) \le 2 \exp(-t^2/B^2)$ for all $t > 0$.
The minimal $B$ is called the {\em subgaussian moment} of $\xi$.
The class of subgaussian random variables includes, among others,
normal, symmetric $\pm 1$, and in general all bounded random variables.

\begin{theorem}[Subgaussian]   \label{t: subgaussian}
  Let $A$ be an $n \times n$ matrix whose entries are i.i.d. centered
  random variables with unit variance and subgaussian moment bounded by $B$.
  Then for every $K \ge 2$ one has
  \begin{equation}                      \label{eq main}
    \P \big( s_n(A) > K n^{-1/2} \big)
    \le (C/K) \log K +  c^n,
  \end{equation}
  where $C > 0$  and $c \in (0,1)$ depend (polynomially) only on $B$.
\end{theorem}

\noindent {\bf Remark.}
  A reverse result was proved in \cite{RV square}: for every $\e \ge 0$, one has
  $\P \big( s_n(A) \le \e n^{-1/2} \big) \le C\e +  c^n$.

\medskip

Our argument is an application of the small ball probability bounds
and the structure theory developed in \cite{RV square} and \cite{RV rectangular}.
We shall give a complete proof of Theorem~\ref{t: subgaussian} only;
we leave to the interested reader to modify the argument as in \cite{RV square}
to obtain Theorem~\ref{t: 4}.

\section{Proof of Theorem~\ref{t: subgaussian}}

By $(e_k)_{k=1}^n$ we denote the canonical basis of the Euclidean space $\R^n$
equipped with the canonical inner product $\< \cdot , \cdot \> $
and Euclidean norm $\|\cdot\|_2$.
By $C, C_1, c, c_1, \ldots$ we shall denote positive constants that
may possibly depend only on the subgaussian moment $B$.

Consider vectors $(X_k)_{k=1}^n$ and $(X^*_k)_{k=1}^n$ an
$n$-dimensional Hilbert space $H$.
Recall that the system $(X_k, X_k^*)_{k=1}^n$ is called a
{\em biorthogonal system} in $H$ if $\< X_j^*, X_k\> = \d_{j,k}$
for all $j,k = 1,\ldots, n$. The system is called {\em complete} if
$\Span(X_k) = H$.
The following notation will be used throughout the paper:
\begin{equation}                                        \label{Hk}
  H_k := \Span(X_i)_{i \ne k}, \qquad
  H_{j,k} := \Span(X_i)_{i \not\in \{j,k\}}, \qquad j,k = 1,\ldots,n.
\end{equation}

The next proposition summarizes some elementary and known properties
of biorthogonal systems.

\begin{proposition}[Biorthogonal systems]                      \label{biorthogonal}
  1. Let $A$ be an $n \times n$ invertible matrix with columns $X_k = A e_k$,
  $k=1,\ldots,n$. Define $X_k^* = (A^{-1})^* e_k$.
  Then $(X_k, X_k^*)_{k=1}^n$ is a complete biorthogonal system in $\R^n$.

  2. Let $(X_k)_{k=1}^n$ be a linearly independent system in an $n$-dimensional
  Hilbert space $H$. Then there exist unique vectors $(X_k^*)_{k=1}^n$ such that
  $(X_k, X_k^*)_{k=1}^n$ is a biorthogonal system in $H$. This system is complete.

  3. Let $(X_k, X_k^*)_{k=1}^n$ be a complete biorthogonal system in a Hilbert space $H$.
  Then $\|X_k^*\|_2 = 1/\dist(X_k,H_k)$ for $k=1,\ldots,n$. \qed
\end{proposition}

Without loss of generality, we can assume that $n \ge 2$ and that
$A$ is a.s. invertible (by adding independent normal random variables
with small variance to all entries of $A$).

Let $u,v > 0$. By \eqref{sn}, the following implication holds:
\begin{equation}                    \label{x Ax goal}
 \exists x \in \R^n : \;  \|x\|_2 \le u, \; \|A^{-1}x\|_2 \ge v n^{1/2}
 \quad \text{implies} \quad
  s_n(A) \le (u/v) n^{-1/2}.
\end{equation}
We will now describe how to find such $x$. Consider the columns
$X_k = A e_k$ of $A$ and the subspaces $H_k$, $H_{j,k}$ defined in \eqref{Hk}.
Let $P_1$ denote the orthogonal projection in $\R^n$ onto $H_1$.
We define the vector
$$
x := X_1 - P_1 X_1.
$$
Define $X_k^* = (A^{-1})^* e_k$. By Proposition~\ref{biorthogonal}
$(X_k, X_k^*)_{k=1}^n$
is a complete biorthogonal system in $\R^n$, so
\begin{equation}                            \label{ker P1}
  \ker(P_1) = \Span(X_1^*).
\end{equation}

Clearly, $\|x\|_2 = \dist(X_1,H_1)$.
Conditioning on $H_1$ and using a standard concentration bound,
we obtain
\begin{equation}                            \label{norm x}
  \P(\|x\|_2 > u) \le C e^{-cu^2}, \qquad u > 0.
\end{equation}
This settles the first bound in \eqref{x Ax goal} with high probability.

To address the second bound in \eqref{x Ax goal}, we write
$A^{-1} x = A^{-1} X_1 - A^{-1} P_1 X_1
= e_1 - A^{-1} P_1 X_1$.
Since $P_1 X_1 \in H_1$, the vector $A^{-1} P_1 X_1$ is supported in
$\{2,\ldots,n\}$ and hence is orthogonal to $e_1$. Therefore
\begin{align*}
\|A^{-1} x\|_2^2
  &> \|A^{-1} P_1 X_1\|_2^2
  = \sum_{k=1}^n \< A^{-1} P_1 X_1, e_k \> ^2 \\
  &= \sum_{k=1}^n \< P_1 (A^{-1})^* e_k, X_1 \> ^2
  = \sum_{k=1}^n \< P_1 X_k^*, X_1 \> ^2.
\end{align*}
The first term of the last sum is zero since $P_1 X_1^* = 0$ by \eqref{ker P1}.
We have proved that
\begin{equation}                                    \label{Yk*}
  \|A^{-1} x\|_2^2 \ge \sum_{k=2}^n \< Y_k^*, X_1 \> ^2,
  \quad \text{where} \quad
  Y_k^* := P_1 X_k^* \in H_1, \qquad k=2,\ldots,n.
\end{equation}

\begin{lemma}                           \label{YkXk}
  $(Y_k^*, X_k)_{k=2}^n$ is a complete biorthogonal system in $H_1$.
\end{lemma}

\noindent {\bf Proof.}
By \eqref{Yk*} and \eqref{ker P1},
$Y_k^* - X_k^* \in \ker(P_1) =  \Span(X_1^*)$,
so $Y_k^* = X_k^* - \l_k X_1^*$ for some $\l_k \in \R$ and all $k=2,\ldots,n$.
By the orthogonality of $X_1^*$ to all of $X_k$, $k=2,\ldots,n$, we have
$\< Y_j^*, X_k \> = \< X_j^*, X_k \> = \d_{j,k}$ for all $j,k=2,\ldots,n$.
The biorthogonality is proved.
The completeness follows since $\dim(H_1)=n-1$.
\qed

In view of the uniqueness in Part 2 of Proposition~\ref{biorthogonal},
Lemma~\ref{YkXk} has the following crucial consequence.

\begin{corollary}                                \label{independence a postiori}
  The system of vectors $(Y_k^*)_{k=2}^n$
  is uniquely determined by the system $(X_k)_{k=2}^n$.
  In particular, the system $(Y_k^*)_{k=2}^n$ and the vector $X_1$ are
  statistically independent. \qed
\end{corollary}

By Part 3 of Proposition~\ref{biorthogonal},
$\|Y_k^*\|_2 = 1/\dist(X_k, H_{1,k})$.
We have therefore proved that
\begin{equation}                                \label{norm sum}
  \|A^{-1} x\|_2^2 \ge \sum_{k=2}^n (a_k/b_k)^2,
  \text{ where }
  a_k = \big| \big\langle \frac{Y_k^*}{\|Y_k^*\|_2}, X_1 \big\rangle \big|, \;
  b_k = \dist(X_k, H_{1,k}).
\end{equation}


We will now need to bound $a_k$ above and $b_k$ below.
Without loss of generality, we will do this for $k=2$.

We are going to use a result of \cite{RV rectangular} that states
that random subspaces have no additive structure.
The amount of structure is formalized
by the concept of the least common denominator.
Given parameters $\a > 0$ and $\g \in (0,1)$, the {\em least common denominator}
of a vector $a \in \R^n$ is defined as
$$
\LCD_{\a,\g}(a)
:= \inf \big\{ \theta > 0: \; \dist (\theta a, \Z^N) < \min(\g\|\theta a\|_2,\a) \big\}.
$$
The least common denominator of a subspace $H$ in $\R^n$ is then defined as
$$
\LCD_{\a,\g}(H)
  = \inf \{ \LCD_{\a,\g} (a) :\; a \in H, \|a\|_2 = 1 \}.
$$
Since $H_{1,2}$ is the span of $n-2$ random vectors with i.i.d. coordinates,
Theorem~4.3 of \cite{RV rectangular} yields that
$$
\P \big\{ \LCD_{\a,c} ((H_{1,2})^\perp) \ge e^{cn} \big\} \ge 1 - e^{-cn}
$$
where $\a = c \sqrt{n}$, and $c>0$ is some constant that may only depend
on the subgaussian moment $B$.

On the other hand, note that the random vector $X_2$ is statistically independent
of the subspace $H_{1,2}$. So, conditioning on $H_{1,2}$ and using the standard
concentration inequality, we obtain
$$
\P \big( b_2 = \dist(X_2, H_{1,2}) \ge t \big) \le C e^{-ct^2},
\qquad t > 0.
$$
Therefore, the event
\begin{equation}                                    \label{EE}
  \EE := \big\{ \LCD_{\a,c}((H_{1,2})^\perp) \ge e^{cn}, \; b_2 < t \big\}
  \text{ satisfies } 
  \P(\EE) \ge 1 - e^{-cn} - C e^{-ct^2}.
\end{equation}

Note that the event $\EE$ depends only on $(X_j)_{j=2}^n$.
So let us fix a realization of $(X_j)_{j=2}^n$ for which $\EE$ holds.
By Corollary~\ref{independence a postiori}, the vector $Y_2^*$ is now fixed.
By Lemma~\ref{YkXk}, $Y_2^*$ is orthogonal to $(X_j)_{j=3}^n$. Therefore
$Y^* := Y_2^*/\|Y_2^*\|_2 \in (H_{1,2})^\perp$, and because event $\EE$ holds, we have
$$
\LCD_{\a,c}(Y^*) \ge e^{cn}.
$$

Let us write in coordinates
$a_2  = |\< Y^*, X_1 \> | = |\sum_{i=1}^n Y^*(i) X_1(i)|$ and recall that
$Y^*(i)$ are fixed coefficients with $\sum_{i=1}^n Y^*(i)^2 = 1$, and $X_1(i)$ are
i.i.d. random variables.
We can now apply Small Ball Probability Theorem~3.3 of \cite{RV rectangular}
(in dimension $m=1$) for this random sum. It yields
\begin{equation}                                    \label{prob a small}
  \P_{X_1} (a_2 \le \e)
  \le C(\e + 1/\LCD_{\a,c}(Y^*) + e^{-c_1 n})
  \le C(\e + e^{-c_2 n}).
\end{equation}
Here the subscript in $\P_{X_1}$ means that we the probability is with
respect to the random variable $X_1$ while the other random variables
$(X_j)_{j=2}^n$ are fixed; we will use similar notations later.

Now we unfix all random vectors, i.e. work with $\P = \P_{X_1,\ldots,X_n}$.
We have
\begin{align*}
\P(a_2 \le \e \text{ or } b_2 \ge t)
 &= \E_{X_2,\ldots,X_n} \P_{X_1} (a_2 \le \e \text{ or } b_2 \ge t) \\
 &\le \E_{X_2,\ldots,X_n} \one_\EE \P_{X_1} (a_2 \le \e)
   + \P_{X_2,\ldots,X_n} (\EE^c)
\end{align*}
because $b_2 < t$ on $\EE$. By \eqref{prob a small} and \eqref{EE},
we continue as
\begin{align*}
\P(a_2 \le \e \text{ or } b_2 \ge t)
  &\le C(\e + e^{-c_2 n}) + (e^{-cn} + Ce^{-ct^2}) \\
  &= C_1(\e + e^{-c_3 t^2} + e^{-cn}) := p(\e,t,n).
\end{align*}
Repeating the above argument for any $k \in \{2,\ldots,n\}$
instead of $k=2$, we conclude that
\begin{equation}                            \label{ak/bk}
  \P \big( a_k/b_k \le \e/t \big) \le p(\e,t,n)
  \qquad \text{for } \e > 0, \; t > 0, \; k=2,\ldots,n.
\end{equation}
From this we can easily deduce the lower bound on the sum of $(a_k/b_k)^2$,
which we need for \eqref{norm sum}. This can be done using the
following elementary observation proved by applying  Markov's inequality twice.

\begin{proposition}                         \label{prob sum}
 Let $Z_k \ge 0$, $k=1,\ldots,n$, be random variables.
 Then, for every $\e >0$, we have
 $$
 \P \Big( \frac{1}{n} \sum_{k=1}^n Z_k \le \e \Big)
 \le \frac{2}{n} \sum_{k=1}^n \P(Z_k \le 2\e). \qquad \qed
 $$
\end{proposition}

We use Proposition~\ref{prob sum} for $Z_k = (a_k/b_k)^2$,
along with the bounds \eqref{ak/bk}. In view of \eqref{norm sum}, we obtain
\begin{equation}                            \label{norm A inverse x}
  \P \big( \|A^{-1}x\|_2 \le (\e/t) n^{1/2} \big)
  \le 2 p(4\e,t,n).
\end{equation}
Estimates \eqref{norm x} and \eqref{norm A inverse x} settle the desired bounds
in \eqref{x Ax goal}, and therefore we conclude that
\begin{align*}
\P \big( s_n(A) \le (ut/\e) n^{-1/2} \big)
  &\ge \P \big( \|x\|_2 \le u, \|A^{-1} x\|_2 \ge (\e/t) n^{1/2} \big)\\
  &\ge 1 - C e^{-cu^2} - 2 p(4\e,t,n).
\end{align*}
This estimate is valid for all $\e,u,t > 0$. Choosing
$\e = 1/K$, $u=t=\sqrt{\log K}$, the proof of Theorem~\ref{t: subgaussian} is complete.
\qed


\begin{thebibliography}{00}

\bibitem{E 88} A. Edelman,
  {\em Eigenvalues and condition numbers of random matrices},
  SIAM J. Matrix Anal. Appl. 9 (1988) 543--560

\bibitem{RV square} M. Rudelson, R. Vershynin,
  {\em The Littlewood-Offord Problem and invertibility of random matrices},
  Advances in Mathematics 218 (2008) 600--633

\bibitem{RV rectangular} M. Rudelson, R. Vershynin,
  {\em The smallest singular value of a random rectangular matrix},
  submitted

\bibitem{vN} J. von Neumann,
  {\em Collected works.
  Vol. V: Design of computers, theory of automata and numerical analysis}.
  General editor: A. H. Taub. A Pergamon Press Book The Macmillan Co.,
  New York, 1963

\end{thebibliography}
\end{document}